\documentclass{amsart}

\usepackage{amsthm,amsfonts,amsmath,amssymb}
\usepackage{xcolor}
\usepackage{graphicx,float,enumerate,subfigure,tikz}
\usepackage{url}
\usepackage{mathtools}
\usepackage{enumitem}
\setlist[enumerate,1]{label=(\roman*),font=\normalfont}
\usepackage{hyperref}
\hypersetup{colorlinks=true, citecolor=blue, linkcolor=blue}
\usepackage[noabbrev,capitalise]{cleveref}
 
\newtheorem{theorem}{Theorem}
\newtheorem{corollary}[theorem]{Corollary}
\newtheorem{lemma}[theorem]{Lemma}
\newtheorem{proposition}[theorem]{Proposition}

\crefname{problem}{Problem}{Problems}

\newtheorem*{theorem*}{Theorem}
\newtheorem*{corollary*}{Corollary}

\theoremstyle{definition}

\newtheorem*{claim*}{Claim}
\crefname{claim}{Claim}{Claims}

\theoremstyle{remark}
\newtheorem{remark}[theorem]{Remark}
\crefname{remark}{Remark}{Remarks}

\crefname{example}{Example}{Examples}

\newcommand{\card}[1]{\ensuremath{|#1|}}

\newcommand{\R}{\mathbb{R}}

\newcommand{\E}{\mathcal{E}}

\newcommand{\ssm}{\smallsetminus}

\newcommand{\defn}[1]{\emph{\color{blue} #1}}
\newcommand{\eqdef}{\mbox{\,\raisebox{0.2ex}{\scriptsize\ensuremath{\mathrm:}}\ensuremath{=}\,}} 
\newcommand{\defeq}{\mbox{\,\ensuremath{=}\raisebox{0.2ex}{\scriptsize\ensuremath{\mathrm:}}\,}} 
\newcommand{\set}[2]{\ensuremath{\left\{#1\,\middle|\,#2\right\}}} 

\title{Edge connectivity of simplicial polytopes}

\author{Vincent Pilaud} 
\address[VP]{CNRS \& LIX, École polytechnique}
\email{vincent.pilaud@lix.polytechnique.fr}
\thanks{VP was partially supported by the French projects CAPPS (ANR~17\,CE40\,0018), and CHARMS (ANR~19\,CE40\,0017), and by the French\,--\,Austrian project PAGCAP (ANR~21\,CE48\,0020 \& FWF I 5788).}

\author{Guillermo Pineda-Villavicencio}
\address[GPV]{Federation University Australia, School of Information Technology, Deakin University}
\email{work@guillermo.com.au}

\author{Julien Ugon}
\address[JU]{Federation University Australia, School of Information Technology, Deakin University}
\email{julien.ugon@deakin.edu.au}
\thanks{JU was partially supported by ARC discovery project DP180100602.}


\begin{document}

\begin{abstract} 
We show that the graph of a simplicial polytope of dimension~${d \ge 3}$ has no nontrivial minimum edge cut with fewer than $d(d+1)/2$ edges, hence the graph is $\min\{\delta,d(d+1)/2\}$-edge-connected where $\delta$ denotes the minimum degree.
When~$d = 3$, this implies that every minimum edge cut in a plane triangulation is trivial.
When~$d \ge 4$, we construct a simplicial $d$-polytope whose graph has a nontrivial minimum edge cut of cardinality $d(d+1)/2$, proving that the aforementioned result is best possible.
\end{abstract}

\maketitle 

\section{Introduction}

A \defn{polytope} is the convex hull of finitely many points, its \defn{dimension} is the dimension of its affine hull, its \defn{faces} are its intersections with its supporting hyperplanes (and the polytope itself), and its \defn{graph} is the graph whose vertices are its $0$-dimensional faces and whose edges are its $1$-dimensional faces. A polytope of dimension $d$ is referred to as a \defn{$d$-polytope}.
A polytope is \defn{simplicial} (resp.~\defn{cubical}) when all its faces are simplices (resp.~combinatorially equivalent to cubes).
See~\cite{Gru03,Zie95} for standard textbooks on polytope theory.

The \defn{vertex} (resp.~\defn{edge}) \defn{connectivity} of a connected graph~$G$ is the minimum number of vertices (resp.~edges) whose removal disconnects $G$.
More precisely, the edge connectivity of~$G$ is the cardinality of the smallest edge cut of~$G$.
An \defn{edge cut} is the set of edges~$\E(V_\circ, V_\bullet)$ with a vertex~in~$V_\circ$ and a vertex~in~$V_\bullet$ for some partition~$V_\circ \sqcup V_\bullet = V$ of the vertex set~$V$ of~$G$ with~${V_\circ \ne \varnothing \ne V_\bullet}$.
A \defn{minimum edge cut} is an edge cut of minimal cardinality.
See~\cite{BM08} for a textbook on graph~theory.

A famous result of Balinski~\cite{Bal61} ensures that the graph of a $d$-polytope is \mbox{$d$-vertex-connected}, which implies that it is $d$-edge-connected. This is the best possible lower bound for the edge connectivity of general polytopes (think about the prism over a simplex).
In contrast, the last two authors have investigated further vertex and edge connectivity properties of graphs of cubical polytopes~\cite{ThiPinUgo20, BuiPinUgo20}.
In this paper, we focus on graphs of simplicial polytopes and establish the \mbox{following property}.

\begin{theorem*}
If $V_\circ \sqcup V_\bullet$ partitions the vertices of a simplicial $d$-polytope and~${\card{V_\bullet} \! \ge \! d \! \ge \! 3}$, the  edge cut~$\E(V_\circ, V_\bullet)$ has at least $m(2d+1-m)/2$ edges where $m = \min(d, \card{V_\circ})$.
\end{theorem*}

Our proof relies on the famous Lower Bound Theorem of Barnette~\cite{Bar71, Bar73} which states (in particular) that the classical stacked polytopes have the minimum number of edges among all simplicial polytopes with the same dimension and number of vertices.
We also use these stacked polytopes to construct simplicial $d$-polytopes that prove that our theorem is best possible for simplicial polytopes.

This theorem implies surprising results on the edge connectivity of the graphs of simplicial polytopes.
Note that the edges incident to a given vertex always form an edge cut that we call \defn{trivial}.

\begin{corollary*}
Every nontrivial minimum edge cut in the graph of a simplicial polytope of dimension~$d \ge 3$ has at least $d(d+1)/2$ edges.
\end{corollary*}

\begin{corollary*}
The graph of a simplicial polytope of dimension~$d \ge 3$ and minimum degree~$\delta$ is $\min\{\delta,d(d+1)/2\}$-edge-connected.
\end{corollary*}

Consider now the specific case of~$d = 3$.
By Steinitz's theorem \cite{Ste22}, the graph of a $3$-polytope is planar and $3$-connected, and a planar realization of the graph of a simplicial 3-polytope is a \defn{plane triangulation}.
Therefore, abusing terminology slightly, we will use interchangeably the terms plane triangulation and simplicial 3-polytope.
In this case, Euler's formula \cite{Eul58,Eul58a} implies that the minimum degree of the plane triangulation is at most five, from which we derive the following statement.
 
\begin{corollary*}
Every minimum edge cut in a plane triangulation is trivial.
\end{corollary*}

As this result is the base case of the inductive proof of our theorem, we provide a short, graph-theoretical proof of this corollary in \cref{sec:planeEdgeConnectivity}.
Surprisingly, this result seems to be new.

\bigskip
We conclude this introduction with a quick historical remark on the genesis of the current version of the paper. In a first version~\cite{PinUgo21}, GPV and JU proved that every nontrivial minimum edge cut in the graph of a simplicial $d$-polytope has at least $4d-8$ edges. The proof relied on links of vertices in a simplicial polytope, but could not be pushed to get a quadratic bound. They also already constructed the examples of simplicial $d$-polytopes with nontrivial edge cuts of cardinality~$d(d+1)/2$. While anonymously refereeing the paper, VP proposed the use of the Lower Bound Theorem for simplicial polytopes \cite{Bar71, Bar73} to establish the quadratic bound $d(d+1)/2$. In agreement with the editors, we decided to combine the results to produce a new version of the paper. Another referee then catched an embarrassing flaw in the proof, forcing us to seriously revise the argument, and leading to the current version of the paper.

\section{Plane triangulations} 
\label{sec:planeEdgeConnectivity}

We first provide elementary proofs of the aforementioned statements in the case of plane triangulations.
Recall that if a simple plane graph~$G$ has~$v \ge 3$ vertices and $e$ edges, then Euler's formula \cite{Eul58,Eul58a} and double counting the edge--face incidences show that $e \le 3 v - 6$, with equality if and only if $G$ is a plane triangulation. 
In turn, this implies that the minimum degree~$\delta$ of~$G$ is at most~$5$, since double counting the vertex--edge incidencies gives~$\delta v \le 2e \le 6v - 12$.

\begin{proposition}
\label{prop:mainDimension3}
If $V_\circ \sqcup V_\bullet$ partitions the vertices of a plane triangulation with ${\card{V_\circ} \ge 1}$ and~${\card{V_\bullet} \! \ge \! 3}$, then the edge cut~$\E(V_\circ, V_\bullet)$ has at least $m(7-m)/2$ edges where $m = \min(3, \card{V_\circ})$.
\end{proposition}

\begin{proof}
Denote by~$G$ the plane triangulation, and by~$G_\circ$ and~$G_\bullet$ the subgraphs of~$G$ induced by~$V_\circ$ and~$V_\bullet$ respectively.
As subgraphs of planar graphs, both $G_\circ$ and $G_\bullet$ are planar.
Denote by~$v, v_\circ, v_\bullet$ the number of vertices and by $e, e_\circ, e_\bullet$ the number of edges of~$G, G_\circ, G_\bullet$, and by~$e_{\circ\bullet}\eqdef \card{\E(V_\circ, V_\bullet)}$ the number of edges of~$G$ from~$V_\circ$ to~$V_\bullet$.
Since~$v_\bullet \ge 3$, we have~${e_\bullet \le 3 v_\bullet - 6}$.
Since~$v = v_\circ + v_\bullet$ and~$e = e_\circ + e_\bullet + e_{\circ\bullet}$, we get
\[
e_{\circ\bullet} = e - e_\circ - e_\bullet \ge (3 v - 6) - e_\circ - (3 v_\bullet - 6) = 3 v_\circ - e_\circ.
\]
Hence,
\begin{itemize}
\item if~$v_\circ = 1$, then~$e_\circ = 0$ and $e_{\circ\bullet} \ge 3 = 1(7-1)/2$,
\item if~$v_\circ = 2$, then~$e_\circ \le 1$ and $e_{\circ\bullet} \ge 5 = 2(7-2)/2$,
\item if~$v_\circ \ge 3$, then~$e_\circ \le 3 v_\circ - 6$ so that~$e_{\circ\bullet} \ge 6 = 3(7-3)/2$. \qedhere
\end{itemize}
\end{proof}

\begin{corollary}
Every minimum edge cut in a plane triangulation is trivial.
\end{corollary}

\begin{proof}
Denoting by~$\delta$ the minimum degree of the triangulation (hence ${3 \le \delta \le 5}$), and using the notations from the previous proof,
\begin{itemize}
\item if~$v_\circ = 1$, then the cut is trivial,
\item if~$v_\circ = 2$, then~$e_\circ \le 1$ and $e_{\circ\bullet} \ge 2 \delta - 1 > \delta$, so that the cut is not minimum,
\item if~$v_\circ \ge 3$, then~$e_{\circ\bullet} \ge 6 > \delta$, so that the cut is not minimum. \qedhere
\end{itemize}
\end{proof}

\section{Simplicial polytopes}
\label{sec:simplicialEdgeConnectivity}

We now consider a partition~$V_\circ \sqcup V_\bullet$ of the vertices of a simplicial $d$-polytope~$P$ with~$V_\circ \ne \varnothing \ne V_\bullet$.
We use the same notations as before: 
\begin{itemize}
\item for vertices, $v = v_\circ + v_\bullet$ where~$v_\circ \eqdef \card{V_\circ}$ and~$v_\bullet \eqdef \card{V_\bullet}$, and $m \eqdef \min(d, v_\circ)$,
\item for edges~$e = e_\circ + e_\bullet + e_{\circ\bullet}$ where~$e_\circ \eqdef \card{\E(V_\circ, V_\circ)}$, $e_\bullet \eqdef \card{\E(V_\bullet, V_\bullet)}$, and $e_{\circ\bullet} \eqdef \card{\E(V_\circ, V_\bullet)}$.
\end{itemize}
Our objective is the following statement, announced in the introduction.

\begin{theorem}
\label{thm:mainPrecise}
If~$v_\bullet \ge d \ge 3$, then $e_{\circ\bullet} \ge m(2d+1-m)/2$.
\end{theorem}

\begin{remark}
Note that if~$v_\circ = 1$, then~$e_{\circ\bullet}$ is the degree of the only vertex of~$V_\circ$, which is indeed at least~$d$.
When~$v_\circ = 2$, \cref{thm:mainPrecise} is equivalent to the fact that if two adjacent vertices of degree~$d$ in a $d$-polytope are contained only in simplex faces, then the polytope is a polygon or a $d$-simplex.
This is a slight refinement of the classical exercise asserting that the polygons and the simplices are the only simple and simplicial polytopes.
We omit the elementary proof here as we do not need it to show \cref{thm:mainPrecise}.
\end{remark}

\subsection{Basic case from the Lower Bound Theorem}
\label{subsec:LBT}

We first observe that \cref{thm:mainPrecise} holds when both parts of the partition are small enough.
The proof relies on the classical Lower Bound Theorem for edges of simplicial polytopes by Barnette \cite{Bar71, Bar73}, which will also be used in the inductive proof of the general case.

\begin{theorem}[\cite{Bar73}] 
\label{thm:LBT} 
A simplicial $d$-polytope with $v$ vertices has at least $dv-\binom{d+1}{2}$ edges.
\end{theorem}

\begin{proposition}
\label{prop:specialCase}
\cref{thm:mainPrecise} holds when~$\max(v_\circ, v_\bullet) \le d+1$.
\end{proposition}

\begin{proof}
Since there are at most $\binom{v_\circ}{2}$ edges inside $V_\circ$ and at most $\binom{v_\bullet}{2}$ inside $V_\bullet$, an application of the Lower Bound Theorem of \cref{thm:LBT} yields that
\begin{align*}
e_{\circ\bullet} & = e - e_\circ - e_\bullet \ge d (v_\circ + v_\bullet) - {\textstyle \binom{d + 1}{2} - \binom{v_\circ}{2} - \binom{v_\bullet}{2}} \\ & = v_\circ(2d+1-v_\circ)/2 + v_\bullet(2d+1-v_\bullet)/2 - d(d+1)/2 \\ & = v_\circ(2d+1-v_\circ)/2 = m(2d+1-m)/2
\end{align*}
where the penultimate equality holds since~$v_\bullet = d$ or~$v_\bullet = d+1$, and the last equality holds since~$v_\circ \le d+1$.
\end{proof}

\subsection{Three auxiliary lemmas}
\label{subsec:preliminaries}

We now prove three auxiliary lemmas that will be essential in our inductive proof of \cref{thm:mainPrecise}.

Consider a vertex~$w$ of~$V_\bullet$ with $\bar{v}_\circ$ neighbors in~$V_\circ$ and~$\bar{v}_\bullet$ neighbors in~$V_\bullet$.
We say that~$w$ is \defn{miscolored} if~$m-\bar{v}_\circ < d-\bar{v}_\bullet$.

\begin{lemma}
\label{lem:miscolored}
\cref{thm:mainPrecise} holds as soon as it holds when no vertex of~$V_\bullet$ is miscolored.
\end{lemma}

\begin{proof}
Fix a simplicial $d$-polytope, and for a partition~$V_\circ \sqcup V_\bullet$ with~$v_\bullet \ge d$ of its vertices define~$k_{\circ\bullet}: = {e_{\circ\bullet} - m(2d+1-m)/2}$. \Cref{thm:mainPrecise} is equivalent to $k_{\circ\bullet}$ being nonnegative for all possible partitions, and thus we can assume that $k_{\circ\bullet}$ is minimum across all possible such partitions.
By \cref{prop:specialCase}, \cref{thm:mainPrecise} holds if~$\max(v_\circ, v_\bullet) \le d+1$.
We can thus assume that~$v_\circ > d$ or~$v_\bullet > d$.
We will now prove that no vertex of~$V_\bullet$ is miscolored, so that \cref{thm:mainPrecise} indeed holds by assumption.

Assume first that~$v_\bullet > d$ and that~$w \in V_\bullet$ is miscolored.
Consider the partition~$V_\circ' \sqcup V_\bullet'$ where~$V_\circ' \eqdef V_\circ \cup \{w\}$ and~$V_\bullet' \eqdef V_\bullet \ssm \{w\}$.
Let~$v_\circ' \eqdef \card{V_\circ'} = v_\circ + 1$, $v_\bullet' \eqdef \card{V_\bullet'} = v_\bullet - 1$, $m' \eqdef \min(d, v_\circ')$ and $e_{\circ\bullet}' \eqdef \card{\E(V_\circ', V_\bullet')}$.
Note that~$v_\bullet' = v_\bullet-1 \ge d$ so that the partition is valid.
As the only modified vertex is~$w$, we have
\[
e_{\circ\bullet} = e_{\circ\bullet}' + \bar{v}_\circ - \bar{v}_\bullet > e_{\circ\bullet}' + m - d,
\]
where the last inequality holds since~$w$ is miscolored. If~$v_\circ < d$ then~$m = v_\circ$ and ${m' = v_\circ' = v_\circ + 1}$, hence~$e_{\circ\bullet} > e_{\circ\bullet}' + v_\circ - d$ so that
\begin{align*}
e_{\circ\bullet} - m(2d+1-m)/2 & > e_{\circ\bullet}' + v_\circ - d - v_\circ(2d+1-v_\circ)/2 \\ & = e_{\circ\bullet}' - (v_\circ + 1)(2d-v_\circ)/2 \\ & = e_{\circ\bullet}' - m'(2d+1-m')/2.
\end{align*}
If~$v_ \circ \ge d$, then~$m = d$ and~$m' = d$, hence~$e_{\circ\bullet} > e_{\circ\bullet}'$ so that
\[
e_{\circ\bullet} - m(2d+1-m)/2 > e_{\circ\bullet}' - m'(2d+1-m')/2.
\]
In both cases, we have contradicted that~$V_\circ \sqcup V_\bullet$ is the partition with minimal~$k_{\circ\bullet}$.

Assume now that~$v_\circ > d$ while~$v_\bullet = d$ and that~$w \in V_\bullet$ is miscolored.
Consider the partition~$V_\circ' \sqcup V_\bullet'$ where~$V_\circ' \eqdef V_\bullet \ssm \{w\}$ and~$V_\bullet' \eqdef V_\circ \cup \{w\}$.
Let~${v_\circ' \eqdef \card{V_\circ'} = v_\bullet - 1}$, $v_\bullet' \eqdef \card{V_\bullet'} = v_\circ + 1$, $m' \eqdef \min(d, v_\circ')$ and $e_{\circ\bullet}' \eqdef \card{\E(V_\circ', V_\bullet')}$.
Note that~$v_\bullet' = v_\circ+1 \ge d$ so that the partition is valid.
As before, we have
\[
e_{\circ\bullet} = e_{\circ\bullet}' + \bar{v}_\circ - \bar{v}_\bullet > e_{\circ\bullet}' + m - d,
\]
where the last inequality holds since~$w$ is miscolored.
Moreover $m = d$ and ${m' = d-1}$ so that $m(2d+1-m) = d(d+1) > (d-1)(d+2) = m'(2d+1-m')$, and we obtain
\[
e_{\circ\bullet} - m(2d+1-m)/2 > e_{\circ\bullet}' - m'(2d+1-m')/2,
\]
contradicting again that~$V_\circ \sqcup V_\bullet$ is the partition with minimal~$k_{\circ\bullet}$.
\end{proof}

\begin{lemma}
\label{lem:neighbors1}
\cref{thm:mainPrecise} holds as soon as it holds when each vertex of~$V_\circ$ is adjacent to at least one vertex of~$V_\bullet$.
\end{lemma}

\begin{proof}
By induction on the number of vertices of~$V_\circ$ with no neighbor in~$V_\bullet$.
Assume that there is such a vertex~$w$.
Observe that~$v_\circ \ge d+1$ (since all  the neighbours of $w$ are in $V_\circ$) and that~$v_\bullet \ge d$ (by assumption in \cref{thm:mainPrecise}).
Let~${V_\circ' \eqdef V_\circ \ssm \{w\}}$ and~$V_\bullet' \eqdef V_\bullet$, and consider the convex hull~$P'$ of~$V_\circ' \sqcup V_\bullet'$.
Since the original polytope~$P$ is simplicial, a slight perturbation of its vertices preserves its graph.
We can thus assume that its vertices are in \defn{general position} (any $d+1$ of its vertices are affinely independent), so that the resulting polytope~$P'$ is still simplicial.
As $\card{V_\circ'} = \card{V_\circ} - 1 \ge d$ and~$\card{V_\bullet'} = \card{V_\bullet} = v_\bullet \ge d$, and there is one less vertex in~$V_\circ'$ with no neighbor in~$V_\bullet'$, we obtain by induction that
\(
e_{\circ\bullet} = \card{\E(V_\circ, V_\bullet)} = \card{\E(V_\circ', V_\bullet')} \ge d(d+1)/2.
\)
\end{proof}

\begin{lemma}
\label{lem:neighbors2}
For any edge~$w_\circ w_\bullet$ with~$w_\circ \in V_\circ$ and~$w_\bullet \in V_\bullet$, the number of neighbors of~$w_\circ$ in~$V_\bullet$ plus the number  of neighbors of~$w_\bullet$ in~$V_\circ$ is at least~$d+1$.
\end{lemma}

\begin{proof}
The edge~$w_\circ w_\bullet$ is contained in at least $d-1$ $2$-faces, which are triangles since the polytope is simplicial.
For each such triangle~$w_\circ w_\bullet w$, the vertex~$w$ is a neighbor of both~$w_\circ$ and~$w_\bullet$ and belongs to either~$V_\circ$ or~$V_\bullet$.
Every such vertex $w$ is counted once when adding the number of neighbors of~$w_\circ$ in~$V_\bullet$ and the number  of neighbors of~$w_\bullet$ in~$V_\circ$. Adding $w_\bullet$ and $w_\circ$ to these $d-1$ vertices $w$, we obtain the desired count.
\end{proof}

\subsection{Inductive proof}
\label{subsec:inductiveProof}

We are now ready to deal with the general case.

\begin{proof}[Proof of \cref{thm:mainPrecise}]
The proof works by induction on~$d \ge 3$.
The base case $d = 3$ was already proved in \cref{prop:mainDimension3}.

We consider a vertex~$w$ of~$V_\bullet$ with $\bar{v}_\circ$ neighbors in~$V_\circ$ and~$\bar{v}_\bullet$ neighbors in~$V_\bullet$.
Moreover, we assume that~$w$ is chosen so that~$\bar{v}_\circ$ is maximal.
By \cref{lem:neighbors1}, we may assume that each vertex~$w_\circ$ in~$V_\circ$ is adjacent to at least one neighbor~$w_\bullet$ in~$V_\bullet$.
Since each ~$w_\bullet\in V_\bullet$ has at most~$\bar{v}_\circ$ neighbors in~$V_\circ$ by our maximality assumption, we obtain that $w_\circ$ has at least~$d+1-\bar{v}_\circ$ neighbors in~$V_\bullet$ by \cref{lem:neighbors2}.

We now consider the \defn{vertex figure}~$\bar{P}$ of~$w$ in~$P$, that is, the polytope obtained by intersecting~$P$ with a hyperplane separating~$w$ from all other vertices of~$P$.
Since~$P$ is a simplicial $d$-polytope, $\bar{P}$ is a simplicial $(d-1)$-polytope, with
\begin{itemize}
\item a vertex~$\bar{x}$ for each neighbor~$x$ of~$w$ in~$P$, and
\item an edge~$\bar{x} \bar{y}$ for each 2-face $xyw$ of~$P$.
\end{itemize}
Consider the partition~$\bar{V}_\circ \sqcup \bar{V}_\bullet$ of the vertices of~$\bar{P}$, where~$\bar{V}_\circ = \set{\bar{x}}{x \in V_\circ \text{ neighbor of } w}$ and~$\bar{V}_\bullet = \set{\bar{x}}{x \in V_\bullet \text{ neighbor of } w}$.
Note that~$\card{\bar{V}_\circ} = \bar{v}_\circ$ and~$\card{\bar{V}_\bullet} = \bar{v}_\bullet$, and denote by~$\bar{e}_{\circ\bullet} \eqdef \card{\E(\bar{V}_\circ, \bar{V}_\bullet)}$ the number of edges of~$\bar{P}$ between~$\bar{V}_\circ$ and~$\bar{V}_\bullet$.

Observe that the following subsets of the cut~$\E(V_\circ, V_\bullet)$ are pairwise disjoint:
\begin{itemize}
\item edges incident to~$w$ ($\bar{v}_\circ$ of them),
\item edges not incident to $w$ that lie in a 2-face of $P$ containing $w$ (these correspond to the $\bar{e}_{\circ\bullet}$ edges of~$\bar{P}$ between $\bar{V}_\circ$ and $\bar{V}_\bullet$),
\item edges incident to some~$w_\circ$ in~$V_\circ$ not adjacent to~$w$ (there are $v_\circ - \bar{v}_\circ$ such vertices, each of which is incident to at least $d+1-\bar{v}_\circ$ edges from the cut~$\E(V_\circ, V_\bullet)$).
\end{itemize}
Hence, we obtain that
\begin{equation}
\label{eq:mainPrecise-1}
  e_{\circ\bullet} \ge \bar{v}_\circ + \bar{e}_{\circ\bullet} + (v_\circ - \bar{v}_\circ)(d+1-\bar{v}_\circ).
\end{equation}

We now distinguish two cases, depending on whether~$\bar{v}_\bullet$ is less than or greater than~${d-1}$.

\medskip
\paragraph{\bf Case 1: $\bar{v}_\bullet \ge d-1$}
In this case, we can directly apply the induction on~$\bar{P}$.
Defining~$\bar{m} \eqdef \min(d-1, \bar{v}_\circ)$, we have
\[
\bar{e}_{\circ\bullet} \ge \bar{m}(2d-1-\bar{m})/2.
\]

If~$\bar{v}_\circ \ge d$, we have~$\bar{m} = d-1$ and we obtain that $\bar{e}_{\circ\bullet}\ge d(d-1)/2$, which together with \eqref{eq:mainPrecise-1} yields that
\(
e_{\circ\bullet} \ge d + d(d-1)/2 = d(d+1)/2
\),
as desired.

If~$\bar{v}_\circ \le d-1$, we have~$\bar{m} = \bar{v}_\circ$, and thus \eqref{eq:mainPrecise-1} becomes
\begin{align*}
e_{\circ\bullet} & \ge \bar{v}_\circ + \bar{v}_\circ(2d-1-\bar{v}_\circ)/2 + (v_\circ - \bar{v}_\circ)(d+1-\bar{v}_\circ) \\ & =  \bar{v}_\circ(\bar{v}_\circ-1)/2 + v_\circ(d+1-\bar{v}_\circ) \\ & \ge \bar{v}_\circ(\bar{v}_\circ-1)/2 + m(d+1-\bar{v}_\circ) \\ & \ge m(2d+1-m)/2.
\end{align*}
To see the last inequality, define
\begin{align*}
f(t) & \eqdef t(t-1)/2 + m(d+1-t) - m(2d+1-m)/2 \\ & = t(t-1)/2 + m(m+1-2t)/2.
\end{align*}
The last inequality amounts to proving that $f(\bar {v}_\circ) \ge 0$.
For this, observe that when~$t \le m$, we have
\(
f'(t) = t - 1/2 - m < 0
\)
hence
\(
f(t) \ge f(m) = 0. 
\)
Since~$\bar{v}_\circ \le v_\circ$ and~$\bar{v}_\circ \le d-1$, we have~$\bar{v}_\circ \le \min(d, v_\circ) = m$, and we conclude that~$f(v_\circ) \ge 0$.

\medskip
\paragraph{\bf Case 2: $\bar{v}_\bullet < d-1$}
In this case, we will need one more careful analysis, since we cannot directly apply induction on~$\bar{P}$.
Define
\begin{equation}
\label{eq:change-variable}
\bar{v}_\circ' \eqdef m - \bar{v}_\circ
\qquad\text{and}\qquad
\bar{v}_\bullet' \eqdef d - \bar{v}_\bullet.    
\end{equation}
By \cref{lem:miscolored}, we can assume that no vertex of $V_\bullet$ is miscolored, and in particular that $w$ is not miscolored, that is~$\bar{v}_\circ' \ge \bar{v}_\bullet'$.
Moreover, since~$\bar{v}_\bullet < d-1$, we obtain that~$\bar{v}_\circ' \ge \bar{v}_\bullet' \ge 0$, which implies that
\begin{equation}
\label{eq:ineq-change-variable}
\bar{v}_\circ'(\bar{v}_\circ'-1) \ge \bar{v}_\bullet'(\bar{v}_\bullet'-1)
\qquad\text{and}\qquad
\bar{v}_\circ \le m \le d.   
\end{equation}
Applying the Lower Bound Theorem stated in \cref{thm:LBT} to~$\bar{P}$, we obtain that
\begin{align}
\bar{e}_{\circ\bullet} &  \ge (d-1)(\bar{v}_\circ + \bar{v}_\bullet) - {\textstyle \binom{d}{2} - \binom{\bar{v}_\circ}{2} - \binom{\bar{v}_\bullet}{2}} \nonumber \\ 
& = (d-1)(m - \bar{v}_\circ' + d - \bar{v}_\bullet') - {\textstyle \binom{d}{2} - \binom{m-\bar{v}_\circ'}{2} - \binom{d-\bar{v}_\bullet'}{2}} & \nonumber \text{(by \eqref{eq:change-variable})} \\ 
& = m(2d-1+m)/2 - \bar{v}_\circ'(2d-2m+\bar{v}_\circ'-1)/2 - \bar{v}_\bullet'( \bar{v}_\bullet'-1)/2 & \nonumber \\ 
& \ge m(2d-1+m)/2 - \bar{v}_\circ'(d-m+\bar{v}_\circ'-1) & \text{(by \eqref{eq:ineq-change-variable})}
\label{eq:edge-cut-link}
\end{align}
Since~$\bar{v}_\circ \le m \le d$ and $m\le {v}_\circ$, combining \eqref{eq:mainPrecise-1}, \eqref{eq:change-variable}, and \eqref{eq:edge-cut-link},  we thus obtain that
\begin{align*}
e_{\circ\bullet} & \ge \bar{v}_\circ + \bar{e}_{\circ\bullet} + (v_\circ - \bar{v}_\circ)(d+1-\bar{v}_\circ) \ge \bar{v}_\circ + \bar{e}_{\circ\bullet} + (m - \bar{v}_\circ)(d+1-\bar{v}_\circ) \\ & \ge m - \bar{v}_\circ' + m(2d-1+m)/2 - \bar{v}_\circ'(d-m+\bar{v}_\circ'-1) + \bar{v}_\circ'(d+1-m+\bar{v}_\circ') \\ & = m(2d+1-m)/2 + \bar{v}_\circ' \ge  m(2d+1-m)/2. \qedhere
\end{align*}
\end{proof}

\subsection{Edge connectivity}
\label{subsec:edgeConnectivity}

In this section, we show that \cref{thm:mainPrecise} gives the announced quadratic bound on the size of nontrivial minimum edge cuts in simplicial $d$-polytopes.
It is based on the following elementary observation.

\begin{lemma}
\label{lem:bigParts}
If a graph has minimum degree~$\delta$ and a nontrivial minimum edge cut $\E(V_\circ, V_\bullet)$, then both~$V_\circ$ and~$V_\bullet$ have cardinality at least~$\delta$.
\end{lemma} 

\begin{proof}
Define as before~$v_\circ \eqdef \card{V_\circ}$, $v_\bullet \eqdef \card{V_\bullet}$, and~$e_{\circ\bullet} \eqdef \card{\E(V_\circ, V_\bullet)}$.
Since the cut is minimum, we have~$e_{\circ\bullet} \le \delta$.
Because each vertex of~$V_\circ$ is adjacent to at least $\delta$ vertices, of which at most $v_\circ-1$ belong to $V_\circ$, we have ${v_\circ(\delta - v_\circ + 1) \le e_{\circ\bullet}}$.
Hence, we obtain that~${v_\circ(\delta - v_\circ + 1) \le \delta}$, which implies that~$v_\circ \ge \delta$ since $1$ and $\delta$ are the two roots of the quadratic polynomial~$x (\delta - x + 1) - \delta = -(x - 1)(x - \delta)$.
By symmetry, we obtain that~$v_\circ \ge \delta$ and~$v_\bullet \ge \delta$.
\end{proof}

\begin{corollary}
\label{coro:nontrivialMinimalEdgeCut}
Every nontrivial minimum edge cut in the graph of a simplicial polytope of dimension~$d \ge 3$ has at least $d(d+1)/2$ edges.
\end{corollary}

\begin{proof}
By \cref{lem:bigParts}, both parts of the cut have size at least the minimum degree of the graph, hence at least~$d$. The result thus directly follows from \cref{thm:mainPrecise}. 
\end{proof}

\begin{corollary}
The graph of a simplicial polytope of dimension~$d \ge 3$ and minimum degree~$\delta$ is $\min\{\delta,d(d+1)/2\}$-edge-connected.
\end{corollary}

\subsection{A construction of nontrivial minimum edge cuts}
\label{subsec:examples}

To conclude, we construct polytopes that show that the bounds of \cref{thm:mainPrecise,coro:nontrivialMinimalEdgeCut}~are~tight.

We need to recall the definitions of two classical families of polytopes.
\begin{enumerate}
\item The \defn{cyclic $d$-polytope} with~$n\ge d+1$ vertices is the convex hull of $n$ arbitrary points on the moment curve~$t \mapsto (t, t^2, t^3, \dots, t^d)$ of~$\R^d$. Note that the cyclic $d$-polytope is simplicial and achieves the maximal number of $i$-faces among all $d$-polytopes with the same number of vertices for any~$i \le d$, as described by the Upper Bound Theorem for polytopes \cite{McM70}. Moreover, its graph is complete for any~$d \ge 4$.

\item A \defn{stacked $d$-polytope} is either a $d$-simplex or a $d$-polytope with $n \ge d+2$ vertices obtained as the convex hull of a stacked $d$-polytope with $n-1$ vertices together with a point located very close to one of its facets (so that it is visible from this facet, but not from the other facets). Note that all stacked $d$-polytopes are simplicial and achieve the minimum number of $i$-faces among all simplicial $d$-polytopes with the same number of vertices for any~$i \le d$, as described by the Lower Bound Theorem for simplicial polytopes \cite{Bar71, Bar73}.
\end{enumerate}

We will use the following very specific family of stacked polytopes.

\begin{lemma}
\label{lem:stackedPolytope}
For any~$m \le d$, there is a stacked $d$-polytope with~$m+d$ vertices partitioned into~$V_\circ \sqcup V_\bullet$ such that~$V_\circ$ forms an $(m-1)$-face~$F_\circ$ and $V_\bullet$ forms a $(d-1)$-face~$F_\bullet$.
Moreover, the corresponding cut $\E(V_\circ, V_\bullet)$ has size~$ m(2d+1-m)/2$.
\end{lemma}

\begin{proof}
Consider the family of stacked $d$-polytopes $S_j$ for $j = 1, \dots, m$ constructed as follows:
\begin{enumerate}
	\item $S_1$ is a $d$-simplex, with vertices labelled $w^\circ_1, w^\bullet_1, \dots, w^\bullet_d$.
	\item For $j = 2, \dots, m$, the polytope $S_j$ is obtained by stacking a vertex~$w^\circ_j$ on the facet of $S_{j-1}$ with vertices $w^\circ_1, \dots, w^\circ_{j-1}, w^\bullet_j, \dots, w^\bullet_d$. Note that the vertices $w^\circ_1, \dots, w^\circ_j, w^\bullet_{j+1}, \dots, w^\bullet_d$ then form a facet of~$S_j$.
\end{enumerate}
Then $S_m$ is a stacked $d$-polytope with $m+d$ vertices, where~$V_\circ \eqdef \{w^\circ_1, \dots, w^\circ_m\}$ forms an $(m-1)$-face~$F_\circ$ while $V_\bullet \eqdef \{w^\bullet_1, \dots, w^\bullet_d\}$ forms a $(d-1)$-face~$F_\bullet$.
Finally, since stacked polytopes achieve the minimum number of edges described by the Lower Bound Theorem for simplicial polytopes, the corresponding cut~$\E(V_\circ, V_\bullet)$ has size
\[
\card{\E(V_\circ, V_\bullet)} \defeq e_{\circ\bullet} = d(m+d)\textstyle{-\binom{d+1}{2}-\binom{m}{2}-\binom{d}{2}} = m(2d+1-m)/2.
\qedhere
\]
\end{proof}

We first use these stacked polytopes to show that the bound of \cref{thm:mainPrecise}~is~tight.

\begin{proposition}
\label{prop:mainOptimal}
For any~$v_\circ \ge 1$ and~$v_\bullet \ge d \ge 3$, there is a partition~$V_\circ \sqcup V_\bullet$ of the vertices of a stacked $d$-polytope with~$\card{V_\circ} = v_\circ$ and~$\card{V_\bullet} = v_\bullet$, and such that~$e_{\circ\bullet} = m(2d+1-m)/2$ where~$m \eqdef \min(d, v_\circ)$.
\end{proposition}

\begin{proof}
Consider the family of stacked $d$-polytopes $S_{j,k}$ for $j = m, \dots, v_\circ$ and ${k = d, \dots, v_\bullet}$ constructed as follows:
\begin{enumerate}
	\item $S_{m,d}$ is the stacked $d$-polytope with~$m+d$ vertices described in \cref{lem:stackedPolytope}.
	\item For $j = m+1, \dots, v_\circ$, the polytope $S_{j,k}$ is obtained by stacking a vertex~$w^\circ_j$ on the facet of $S_{j-1,k}$ with vertices $w^\circ_{j-d}, \dots, w^\circ_{j-1}$. Note that the vertices $w^\circ_{j-d+1}, \dots, w^\circ_j$ then form a facet of~$S_{j,k}$.
	\item For $k = d+1, \dots, v_\bullet$, the polytope $S_{j,k}$ is obtained by stacking a vertex~$w^\bullet_k$ on the facet of $S_{j,k-1}$ with vertices $w^\bullet_{k-d}, \dots, w^\bullet_{k-1}$. Note that the vertices $w^\bullet_{k-d+1}, \dots, w^\bullet_k$ then form a facet of~$S_{j,k}$.
\end{enumerate}
Note that we do nothing in Step~(ii) when~$v_\circ \le d$.
Observe also that the order in which these stackings are performed is not relevant.
The polytope~$S_{v_\circ, v_\bullet}$ is a staked $d$-polytope whose vertices are partitioned by~$V_\circ \sqcup V_\bullet$ where~$V_\circ \eqdef \set{w^\circ_j}{j \in [v_\circ]}$ and~$V_\bullet \eqdef \set{w^\bullet_k}{k \in [v_\bullet]}$.
Moreover, the edges of the cut~$\E(V_\circ, V_\bullet)$ all belong to the original polytope~$S_{m+d}$.
Hence, the cut has size~$e_{\circ\bullet} = m(2d+1-m)/2$.
\end{proof}

\begin{figure}
	\centerline{\includegraphics[scale=.9]{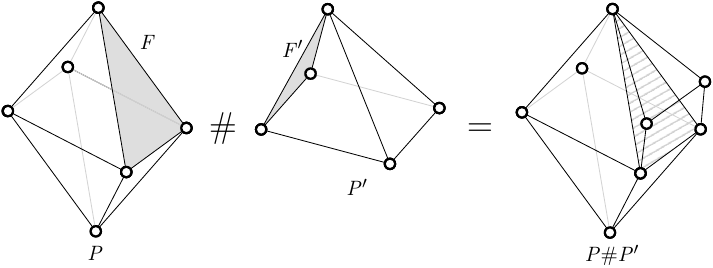}}
	\caption{Connected sum of two polytopes.}
	\label{fig:connected-sum} 
\end{figure}

To prove that the bound of \cref{coro:nontrivialMinimalEdgeCut} is tight, we additionally need the following classical transformations on simplicial polytopes.
Let $P$ and $P'$ be two $d$-polytopes with a facet $F$ of $ P$ projectively isomorphic to a facet $F'$ of $P'$.
Their \defn{connected sum} $P \# P'$ is obtained by ``gluing" $P$ and $P'$ along $F$ and $F'$.
Projective transformations on the polytopes $P$ and $P'$ may be required for the connected sum to be convex.
The connected sum of two polytopes is depicted in \cref{fig:connected-sum}.
This operation was used for instance by Eckhoff to prove that $f$-vectors of polytopes are not unimodal; see \cite[Exm.~8.41]{Zie95}. Observe that stacking a vertex on a simplex facet of a polytope amounts to performing the connected sum of the polytope with a simplex along the facet. 
Our next construction is based on performing connected sums of simplicial polytopes, which is always possible because every polytope combinatorially isomorphic to a simplex is projectively isomorphic to the simplex \cite{McM76}.

\begin{proposition}
\label{prop:nontrivialMinimalEdgeCutOptimal}
For each $d\ge 4$, there is a simplicial $d$-polytope with minimum degree at least $d(d+1)/2$ and a nontrivial minimum edge cut with ${d(d+1)}/{2}$ edges.
\end{proposition}

\begin{proof}
Let~$S$ be a stacked polytope with~$2d$ vertices and two disjoint facets~$F_\circ$ and~$F_\bullet$ (whose existence is guaranteed by \cref{lem:stackedPolytope} when~$m = d$).
Let $C$ be the cyclic $d$-polytope with $1 + d(d+1)/2$ vertices.
Consider the polytope~$P$ obtained by the connected sum~$C \# S \# C$, where the copies of $C$ are glued along the facets~$F_\circ$ and~$F_\bullet$ of~$S$, and along arbitrary facets of the two copies of~$C$.
Since the graph of~$C$ is complete with $1 + d(d+1)/2$ vertices, it has minimum degree~$d(d+1)/2$.
Hence, so does the graph of~$P$.
Moreover, the number of edges joining a vertex of~$F_\circ$ to a vertex of~$F_\bullet$ is precisely~$d(d+1)/2$ by \cref{lem:stackedPolytope}.
We conclude that~$P$ has a nontrivial edge cut of size~$d(d+1)/2$, which is thus minimum by \cref{coro:nontrivialMinimalEdgeCut}.
\end{proof} 
		
\section{Acknowledgements}

We thank Arnau Padrol for various discussions.
We are also grateful to two anonymous referees for various suggestions, in particular for pointing out to us a serious mistake in a previous approach to the proof of \cref{coro:nontrivialMinimalEdgeCut}.


\bibliographystyle{alpha}
\bibliography{simplicial-edge-connectivity}

\end{document}